\documentclass[12pt]{article}
\usepackage{latexsym}
\usepackage{amssymb}
\usepackage{amsmath}
\usepackage{mathrsfs}

\allowdisplaybreaks

\newtheorem{theorem}{Theorem}[section]
\newtheorem{lemma}{Lemma}[section]

\newtheorem{definition}{Definition}[section]
\newtheorem{remark}{Remark}[section]

\begin{document}

\title{\bf Stochastic inertial manifolds for damped wave equations\footnote{This work is
partially supported by the 985 project of Jilin University.}}

\author{Zhenxin Liu\\
{\small College of Mathematics, Jilin University,  Changchun 130012,
People's Republic of China}\\
{\small zxliu@email.jlu.edu.cn}}

\date{\today}
\maketitle

\begin{abstract}
In this paper, stochastic inertial manifold for damped wave
equations subjected to additive white noise is constructed by the
Lyapunov-Perron method. It is proved that when the intensity of
noise tends to zero the stochastic inertial manifold converges to
its deterministic counterpart almost surely.
\end{abstract}

\section{Introduction}

The inertial manifold (IM) introduced by Constantin, Foias,
Nicolaenko, Sell, and Temam \cite{Fo,Foi,FNST0,FNST,Con0} is a
finite dimensional Lipschitz invariant manifold attracting solutions
exponentially, which goes back to the works of Ma\~n\'e, Henry and
Mora \cite{Man,Hen,Mor}. Global attractor is an invariant compact
set attracting solutions which often has a finite (fractal)
dimension and, therefore, it is an important object for the study of
long time behavior of evolution equations. At the present level of
understanding of dynamical systems, global attractors are expected
to be very complicated objects (fractals) and their practical
utilization, for instance for numerical simulations, may be
difficult. The IMs, when they exist, are more convenient objects
which are able to describe the large-time behavior of dynamical
systems. One of the important properties of inertial manifolds is
that they contain global attractors, so the study of dynamics of
infinite dimensional nonlinear systems can be reduced to the study
of dynamics of flows on the inertial manifold, which, in turn, is
described by the dynamics of an ordinary differential equation.
There are extensive works on IMs. See, for example, Chow and Lu
\cite{Cho}, Chow et al \cite{Cho1}, Constantin et al \cite{Con0},
Constantin et al \cite{Con}, Foias et al \cite{FNST0,FNST}, Foias et
al \cite{Fo,Foi}, Foias et al \cite{Foi1}, Mallet-Paret and Sell
\cite{Mal}, Sell and You \cite{Sel}, Temam \cite{Tem}, among others.

Stochastic partial differential equations (SPDE) have been drawing
more and more attention for their importance in describing many
natural phenomenon under random influences. With the rapid
development of random dynamical systems (RDS) \cite{Ar1}, many SPDEs
are studied in the framework of RDS. On many occasions, the
development of SPDE and RDS mimics the deterministic case and many
efforts are devoted to establish the results for SPDE and RDS
corresponding to that for the deterministic case. This is true for
IM: there have been some works on Stochastic IMs, see, for example,
Bensoussan and Flandoli \cite{Ben}, Chueshov and Girya \cite{Chu},
Chueshov and Scheutzow \cite{Chu1}, Da Prato and Debussche
\cite{Dap}, Duan et al \cite{Dua,Dua1}. These works mainly deal with
stochastic parabolic equations. In present paper, we aim to obtain
the existence of stochastic IM for damped wave equations subjected
to additive white noise. Moreover, we will show that the stochastic
IM converges to its deterministic counterpart {\em almost surely}
when the intensity of noise tends to zero. As in the deterministic
case, the usual methods to obtain the existence of stochastic IMs
are Hadamard's graph transform method \cite{Had} and
Lyapunov-Perron's method \cite{Lia,Per}. In this paper, we adopt the
latter one. In Section 2 we introduce some preliminaries and prove
the existence theorem of stochastic IM for abstract evolution
equations with random coefficients; in Section 3 we apply the result
of Section 2 to damped wave equations subjected to additive white
noise and study the property of its IM.

\section{Existence of Stochastic IM for abstract equations}

\begin{definition}
Let $X$ be a metric space with a metric $d_X$. A {\em random
dynamical system (RDS)}, shortly denoted by
$\varphi$, consists of two ingredients: \\
(i) A model of the noise, namely a metric dynamical system $(\Omega,
\mathscr F, \mathbb P, (\theta_t)_{t\in \mathbb R})$, where
$(\Omega, \mathscr F, \mathbb P)$ is a probability space and $(t,
\omega)\mapsto \theta_t\omega$ is a measurable flow which leaves
$\mathbb P$ invariant, i.e. $\theta_t\mathbb P=\mathbb P$ for all
$t\in \mathbb R$.\\
(ii) A model of the system perturbed by noise, namely a cocycle
$\varphi$ over $\theta$, i.e. a measurable mapping $\varphi: \mathbb
R^+\times \Omega\times X \rightarrow X,
(t,\omega,x)\mapsto\varphi(t,\omega,x)$, such that:
\begin{equation}\label{phi}
\varphi(0,\omega,\cdot)={\rm id}_X,
\varphi(t+s,\omega,\cdot)=\varphi(t,\theta_s
\omega,\varphi(s,\omega,\cdot))\quad {\rm for ~all}\quad
t,s\in\mathbb R^+,\omega\in\Omega.
\end{equation}
\end{definition}

Although it is well known that a large class of partial differential
equations with stationary random coefficients and It\^o stochastic
ordinary differential equations generate RDS (for details see
Chapter 1 of \cite{Ar1}), this problem is still unsolved for SPDE
with general noise terms, see \cite{Dua} for the reason. Indeed, the
existence of RDS generated by SPDE has been proved in relatively
narrow generality. In fact, only cases in which the SPDE can be
reduced to a deterministic one with random coefficients can be
treated in the framework of RDS. See, for example,
\cite{Car,Cra,Cra1,Dua1}.

For later use, assume $z$ is an Ornstein-Uhlenbeck process which
satisfies the following equation
\begin{equation}\label{z}
{\rm d}z+\lambda z{\rm d}t=\delta{\rm d}W
\end{equation}
for some $\lambda>0$ and $\delta>0$. The process $z$ has the
following properties, see \cite{Car,Dua} for the proof.

\begin{lemma}\label{zpr}
(i) There exists a $\{\theta_t\}_{t\in\mathbb R}$-invariant set
$\Omega\in\mathcal B(C_0(\mathbb R,\mathbb R))$ of full measure with
sublinear growth:
\[
\lim_{t\rightarrow\pm\infty}\frac{|\omega(t)|}{|t|}=0,~\omega\in\Omega.
\]
(ii) For $\omega\in\Omega$ the random variable
\[
z(\omega)=-\lambda\delta\int_{-\infty}^0e^{\lambda\tau}\omega(\tau){\rm
d}\tau
\]
exists and generates a unique stationary solution of (\ref{z}) given
by
\[
(t,\omega)\rightarrow
z(\theta_t\omega)=-\lambda\delta\int_{-\infty}^0e^{\lambda\tau}\theta_t\omega(\tau){\rm
d}\tau=-\lambda\delta\int_{-\infty}^0e^{\lambda\tau}\omega(\tau+t){\rm
d}\tau+\delta\omega(t).
\]
The map $t\rightarrow z(\theta_t\omega)$ is continuous. \\
(iii) In particular, we have
\[
\lim_{t\rightarrow\pm\infty}\frac{|z(\theta_t\omega)|}{|t|}=0~{\rm
for~}\omega\in\Omega.
\]
(iv) In addition,
\[
\lim_{t\rightarrow\pm\infty}\frac1t\int_0^tz(\theta_\tau\omega){\rm
d}\tau=0~{\rm
for~}\omega\in\Omega,~\lim_{t\rightarrow\pm\infty}\frac1t\int_0^t|z(\theta_\tau\omega)|{\rm
d}\tau=\mathbb E|z|<\infty.
\]
\end{lemma}

Let $H$ be a separable Hilbert space with norm $|\cdot|$ and inner
product $\langle\cdot,\cdot\rangle$. Consider the Stratonovich SPDE
on $H$
\begin{equation}\label{u}
\frac{{\rm d}u}{{\rm d}t}=Au+F(u)+u\dot W,
\end{equation}
where $u\in H$, $W(t)$ is the standard real-valued two-sided Wiener
process and the generalized time-derivative $\dot W$ formally
describes a white noise. Here we assume that $F$ is globally
Lipschitz continuous on $H$ with Lipschitz constant ${\rm Lip}F$.
For the existence and uniqueness theory of (\ref{u}) we can first
write it into its equivalent It\^o equation and then refer to
\cite{Dap1} for details. Under the transformation
$T(\omega,u)=ue^{-z(\omega)}$, (\ref{u}) is conjugated to the
following equation with random coefficients
\begin{equation}\label{main}
\frac{{\rm d}u}{{\rm
d}t}=Au+z(\theta_t\omega)u+G(\theta_t\omega,u),~~u(0)=x\in H,
\end{equation}
where $z$ satisfies
\[
{\rm d}z+z{\rm d}t={\rm d}W,
\]
and $G(\omega,u)=e^{-z(\omega)}F(ue^{z(\omega)})$. It is clear that
${\rm Lip}_uG={\rm Lip}F$.

Assume $A: D(A)\rightarrow H$ is a linear operator which generates a
strongly continuous semigroup $e^{At}$ on $H$, which satisfies the
pseudo exponent dichotomy condition with exponents $0>\alpha>\beta$
and bound $K>0$, i.e. there exists a continuous projection $P$ on
$H$ such that\\
(i) $Pe^{At}=e^{At}P$;\\
(ii) the restriction $e^{At}|_{R(P)}, t\ge 0$, is an isomorphism of
the range $R(P)$ of $P$ onto itself, and we denote $e^{At}$ for
$t<0$ the inverse map;\\
(iii)
\begin{equation}\label{expo}
\left.
        \begin{array}{ll}
          |e^{At}Px|\le Ke^{\alpha t}|x|, & t\le 0, \\
          |e^{At}Qx|\le Ke^{\beta t}|x|, & t\ge 0,
        \end{array}
      \right.
\end{equation}
where $Q=I-P$.

\begin{definition}\label{im}
A random set is called {\em invariant} for RDS $\varphi$ if
\[
\varphi(t,\omega,M(\omega))\subset M(\theta_t\omega), ~{\rm
for~any~} t\ge 0.
\]
If an invariant set $M(\omega)$ can be represented by a Lipschitz or
$C^k$ mapping
\[
h(\cdot,\omega):PH\rightarrow QH
\]
such that
\[
M(\omega)=\{\xi+h(\xi,\omega)|\xi\in PH\},
\]
then we call $M(\omega)$ a {\em Lipschitz or $C^k$ invariant
manifold}. Furthermore, if $PH$ is finite dimensional and
$M(\omega)$ attracts exponentially all the orbits of $\varphi$, then
we call $M(\omega)$ a {\em stochastic inertial manifold of
$\varphi$}.
\end{definition}

\begin{theorem}\label{the} (\cite{Dua1}) If
\begin{equation}\label{gap}
K{\rm
Lip}F\left(\frac{1}{\alpha-\eta}+\frac{1}{\eta-\beta}\right)<1,
\end{equation}
then there exists a Lipschitz invariant manifold for the random
evolutionary Equation (\ref{main}), which is given by
\begin{equation}\label{M}
M(\omega)=\{\xi+h(\xi,\omega)|\xi\in PH\},
\end{equation}
where $h:PH\rightarrow QH$ is a Lipschitz continuous mapping given
by
\begin{equation}\label{h}
h(\xi,\omega)=\int_{-\infty}^0e^{-As+\int_s^0z(\theta_r\omega){\rm
d}r}QG(\theta_s\omega,u(s;\xi,\omega)).
\end{equation}

\end{theorem}

\begin{remark}\rm
It is easy to see that if $F$ is $C^1$, then the stochastic
invariant manifold obtained in Theorem \ref{the} is $C^1$ by Theorem
5.3 of \cite{Dua1}.
\end{remark}

Theorem \ref{the} says that (\ref{main}) has a Lipschitz manifold if
the spectral gap condition (\ref{gap}) holds. To show that the
manifold is an inertial manifold for (\ref{main}), we should verify
that it attracts exponentially all the orbits of $\varphi$. A
stronger reduction property is the exponential tracking property
\cite{Foi1}, also called {\em asymptotical completeness property}
\cite{Rob}: each trajectory of the evolution equation tends
exponentially to a trajectory on the inertial manifold. To be more
specific, we states it as follows:

\begin{definition}\label{asymp}
Let $M(\omega)$ be an invariant manifold for RDS $\varphi$. If for
$\forall x\in H$, there exists an $\bar x\in M(\omega)$ such that
\[
|\varphi(t,\omega,x)-\varphi(t,\omega,\bar x)|\le
c_1e^{-c_2t}|x-\bar x|, ~\forall t\ge 0,
\]
where $c_1>0$ is a constant dependent on $\omega$, $x$ and $\bar x$,
while $c_2$ is a constant independent of these variables, then
$M(\omega)$ is said to have the {\em asymptotic completeness
property}.
\end{definition}

If $M(\omega)$ has the asymptotic completeness property, then the
asymptotic behavior of $\varphi$ on $H$ can be reduced to
$M(\omega)$. Hence the the original infinite dimensional SPDE
problem on $H$ is reduced to a finite dimensional stochastic ODE
problem on $M(\omega)$.

Denote
\[
C_\eta^+:=\{\phi:[0,\infty)\rightarrow H|\phi~{\rm
continuous},~\sup_{t\ge 0}e^{-\eta t-\int_0^tz(\theta_r\omega){\rm
d}r}|\phi(t)|<\infty\},
\]
then $C_\eta^+$ is a Banach space with norm
$|\phi|_{C_\eta^+}:=\sup_{t\ge 0}e^{-\eta
t-\int_0^tz(\theta_r\omega){\rm d}r}|\phi(t)|$.

\begin{theorem}\label{th1} If we have the spectral gap condition
\begin{equation}\label{gap1}
K{\rm
Lip}F\left(\frac{1}{\alpha-\eta}+\frac{1}{\eta-\beta}\right)+K^2{\rm
Lip}h\cdot{\rm Lip}F\frac{1}{\alpha-\eta}<1,
\end{equation}
then the Lipschitz invariant manifold for (\ref{main}) obtained in
Theorem \ref{the} has the asymptotic completeness property.
\end{theorem}
{\bf Proof.} Assume $u,\bar u$ are two solutions of (\ref{main}) and
let $w=\bar u-u$, then $w$ satisfies the following equation:
\begin{equation}\label{w}
\frac{{\rm d}w}{{\rm d}t}=Aw+z(\theta_t\omega)w+\tilde
F(\theta_t\omega,w),
\end{equation}
where
\[
\tilde
F(\theta_t\omega,w):=G(\theta_t\omega,u+w)-G(\theta_t\omega,u).
\]
It is clear that
\begin{equation}\label{tildeF}
\tilde F(\theta_t\omega,0)=0,~{\rm Lip}_w\tilde F={\rm Lip}_uG={\rm
Lip}F.
\end{equation}


First if $w\in C_\eta^+$ is a solution of (\ref{w}), then $w$ can be
expressed by
\begin{align}\label{w1}
w(t)=&e^{At+\int_0^tz(\theta_r\omega){\rm
d}r}Qw(0)+\int_0^te^{A(t-s)+\int_s^tz(\theta_r\omega){\rm
d}r}Q\tilde F(\theta_s\omega,w(s)){\rm
d}s\nonumber\\
&~~+\int_{\infty}^te^{A(t-s)+\int_s^tz(\theta_r\omega){\rm
d}r}P\tilde F(\theta_s\omega,w(s)){\rm d}s.
\end{align}
In fact, since $w$ is a solution of (\ref{w}), we have
\[
w(t)=e^{A(t-t_0)+\int_{t_0}^tz(\theta_r\omega){\rm
d}r}w(t_0)+\int_{t_0}^te^{A(t-s)+\int_s^tz(\theta_r\omega){\rm
d}r}\tilde F(\theta_s\omega,w(s)){\rm d}s.
\]
This implies
\[
Pw(t)=e^{A(t-t_0)+\int_{t_0}^tz(\theta_r\omega){\rm
d}r}Pw(t_0)+\int_{t_0}^te^{A(t-s)+\int_s^tz(\theta_r\omega){\rm
d}r}P\tilde F(\theta_s\omega,w(s)){\rm d}s.
\]
When $t_0>t$, by (\ref{expo}) we have
\begin{align*}
|e^{A(t-t_0)+\int_{t_0}^tz(\theta_r\omega){\rm d}r}Pw(t_0)|&\le
Ke^{\alpha(t-t_0)+\int_{t_0}^tz(\theta_r\omega){\rm d}r}|w(t_0)|\\
&\le Ke^{-(\alpha-\eta)t_0+\alpha t+\int_{0}^tz(\theta_r\omega){\rm
d}r}|w|_{C_\eta^+}.
\end{align*}
By the property of $z(\omega)$ we obtain
\[
e^{A(t-t_0)+\int_{t_0}^tz(\theta_r\omega){\rm d}r}Pw(t_0)\rightarrow
0{\rm ~as~}t_0\rightarrow\infty.
\]
Therefore,
\[
Pw(t)=\int_{\infty}^te^{A(t-s)+\int_s^tz(\theta_r\omega){\rm
d}r}P\tilde F(\theta_s\omega,w(s)){\rm d}s.
\]
Thus (\ref{w1}) holds.

We then show that (\ref{w1}) has solutions on $C_\eta^+$ and $\bar
u(0)=u(0)+w(0)\in M(\omega)$. From \cite{Dua1} we know that the
solution $\bar u$ lies on $M$ if and only if $Q\bar u(0)=h(P\bar
u(0),\omega)$, recalling that $M(\omega)=\{\xi+h(\xi,\omega)|\xi\in
PH\}$. That is
\begin{equation}\label{qw0}
Qw(0)=-Qu(0)+h(Pu(0)+Pw(0),\omega).
\end{equation}
Let
\begin{align*}
\tilde Tw(t)&=e^{At+\int_0^tz(\theta_r\omega){\rm d}r}Qw(0),\\
Tw(t)&=\int_0^te^{A(t-s)+\int_s^tz(\theta_r\omega){\rm d}r}Q\tilde
F(\theta_s\omega,w(s)){\rm d}s
+\int_{\infty}^te^{A(t-s)+\int_s^tz(\theta_r\omega){\rm d}r}P\tilde
F(\theta_s\omega,w(s)){\rm d}s,
\end{align*}
then (\ref{w1}) reads as
\[
w(t)=\tilde Tw(t)+Tw(t).
\]
We assert that $\tilde T$ and $T$ map $C_\eta^+$ to $C_\eta^+$. In
fact,
\begin{align*}
e^{-\eta t-\int_0^tz(\theta_r\omega){\rm d}r}|\tilde Tw(t)|& \le
Ke^{-(\eta-\beta)t}|Qw(0)|\\
&\le K|Qw(0)|\\
&\le^{(\ref{qw0})}K\bigl(|-Qu(0)+h(Pu(0),\omega)|\\
&\qquad\qquad +|h(Pu(0)+Pw(0),\omega)-h(Pu(0),\omega)|\bigr)\\
&\le K\bigl(|-Qu(0)+h(Pu(0),\omega)|+{\rm Lip}h|Pw(0)|\bigr)\\
&\le^{(\ref{w1})} K\biggl(|-Qu(0)+h(Pu(0),\omega)|\\
&\qquad\qquad+{\rm
Lip}h\left|\int_\infty^0e^{-As+\int_s^0z(\theta_r\omega){\rm
d}r}P\tilde
F(\theta_s\omega,w(s)){\rm d}s\right|\biggr)\\
&\le^{(\ref{tildeF})} K\bigl(|-Qu(0)+h(Pu(0),\omega)|+K{\rm
Lip}h\cdot{\rm Lip}F\frac1{\alpha-\eta}|w|_{C_\eta^+}\bigr)
\end{align*}
and
\begin{align*}
e^{-\eta t-\int_0^tz(\theta_r\omega){\rm d}r}|Tw(t)|& \le Ke^{-\eta
t-\int_0^tz(\theta_r\omega){\rm
d}r}\biggl(\int_0^te^{\beta(t-s)+\int_s^tz(\theta_r\omega){\rm
d}r}|\tilde
F(\theta_s\omega,w(s))|{\rm d}s\\
&\qquad
+\left|\int_\infty^te^{\alpha(t-s)+\int_s^tz(\theta_r\omega){\rm
d}r}\tilde F(\theta_s\omega,w(s)){\rm d}s\right|\biggr)\\
&\le^{(\ref{tildeF})} K{\rm Lip}Fe^{-\eta
t-\int_0^tz(\theta_r\omega){\rm
d}r}\biggl(\int_0^te^{\beta(t-s)+\int_s^tz(\theta_r\omega){\rm
d}r}|w(s)|{\rm d}s\\
&\qquad +\int_t^\infty e^{\alpha(t-s)+\int_s^tz(\theta_r\omega){\rm
d}r}|w(s)|{\rm d}s\biggr)\\
&\le  K{\rm Lip}F\biggl(\int_0^te^{-(\eta-\beta)(t-s)}{\rm
d}s+\int_t^\infty e^{(\alpha-\eta)(t-s)}{\rm d}s\biggr)|w|_{C_\eta^+}\\
&\le K{\rm
Lip}F\left(\frac1{\eta-\beta}+\frac1{\alpha-\eta}\right)|w|_{C_\eta^+}.
\end{align*}

Next we show that under the spectral gap condition (\ref{gap1}), the
map $\tilde T+T:C_\eta^+\rightarrow C_\eta^+$ is contractive. To
this end, assume $w,\bar w\in C_\eta^+$, then we have
\begin{align*}
e^{-\eta t-\int_0^tz(\theta_r\omega){\rm d}r}|\tilde Tw(t)-\tilde
T\bar w(t)|& \le^{(\ref{qw0})}
Ke^{-(\eta-\beta)t}|h(Pu(0)+Pw(0),\omega)\\
&\qquad\qquad -h(Pu(0)+P\bar w(0),\omega)|\\
&\le K{\rm Lip}h|Pw(0)-P\bar w(0)|\\
&\le K{\rm
Lip}h\left|\int_\infty^0e^{-As+\int_s^0z(\theta_r\omega){\rm
d}r}P{\rm Lip}\tilde F|w(s)-\bar w(s)|{\rm d}s\right|\\
&\le K^2{\rm Lip}h\cdot{\rm
Lip}F\int_0^\infty e^{-(\alpha-\eta)s}|w-\bar w|_{C_\eta^+}{\rm d}s\\
&\le K^2{\rm Lip}h\cdot{\rm Lip}F\frac1{\alpha-\eta}|w-\bar
w|_{C_\eta^+},
\end{align*}
and
\begin{align*}
&e^{-\eta t-\int_0^tz(\theta_r\omega){\rm d}r}|Tw(t)-T\bar w(t)|\\
\le & e^{-\eta t-\int_0^tz(\theta_r\omega){\rm
d}r}\biggl(\int_0^te^{A(t-s)+\int_s^tz(\theta_r\omega){\rm
d}r}Q|\tilde F(\theta_s\omega,w(s))-\tilde F(\theta_s\omega,\bar
w(s))|{\rm
d}s\\&+\left|\int_\infty^te^{A(t-s)+\int_s^tz(\theta_r\omega){\rm
d}r}P|\tilde F(\theta_s\omega,w(s))-\tilde F(\theta_s\omega,\bar
w(s))|{\rm d}s\right|\biggr)\\
\le & K\biggl(\int_0^te^{\beta(t-s)-\eta
t+\int_s^0z(\theta_r\omega){\rm d}r}{\rm Lip}F|w(s)-\bar w(s)|{\rm
d}s\\
&\quad +\left|\int_\infty^te^{\alpha(t-s)-\eta
t+\int_s^0z(\theta_r\omega){\rm d}r}{\rm Lip}F|w(s)-\bar w(s)|{\rm
d}s\right|\biggr)\\
\le & K{\rm Lip}F\biggl(\int_0^te^{-(\eta-\beta)(t-s)}|w-\bar
w|_{C_\eta^+}{\rm d}s+\int_t^\infty e^{(\alpha-\eta)(t-s)}|w-\bar
w|_{C_\eta^+}{\rm d}s\biggr)\\
\le & K{\rm
Lip}F\left(\frac1{\eta-\beta}+\frac1{\alpha-\eta}\right)|w-\bar
w|_{C_\eta^+}.
\end{align*}
That is
\begin{align*}
&|\tilde Tw-\tilde T\bar w|_{C_\eta^+}\le K^2{\rm Lip}h\cdot{\rm
Lip}F\frac1{\alpha-\eta}|w-\bar w|_{C_\eta^+},\\
&|Tw-T\bar w|_{C_\eta^+}\le K{\rm
Lip}F\left(\frac1{\eta-\beta}+\frac1{\alpha-\eta}\right)|w-\bar
w|_{C_\eta^+}.
\end{align*}
Therefore,
\[
|(\tilde T+T)w-(\tilde T+T)\bar w|_{C_\eta^+}\le \left[K^2{\rm
Lip}h\cdot{\rm Lip}F\frac1{\alpha-\eta}+K{\rm
Lip}F\left(\frac1{\eta-\beta}+\frac1{\alpha-\eta}\right)\right]|w-\bar
w|_{C_\eta^+}.
\]
Then by the spectral gap condition (\ref{gap1}) we obtain that
$\tilde T+T$ has a unique fixed point $w^*$ on $C_\eta^+$, which
satisfies
\[
\bar u(0)=u(0)+w^*(0)\in M(\omega)
\]
as desired. Hence
\begin{align*}
|\bar u(t,\omega,\bar u_0)-u(t,\omega,u_0)|& \le e^{\eta
t+\int_0^tz(\theta_r\omega){\rm d}r}|\bar u_0-u_0|\\
&\le c(\omega)e^{\eta t}|\bar u_0-u_0|, ~t\ge 0
\end{align*}
for some $c(\omega)>0$ by the property of $z(\omega)$. \hfill$\Box$

\section{Stochastic IM for wave equations}

Consider the following wave equation in $[0,\pi]$ perturbed by
additive white noise:
\begin{equation}\label{wa}
\epsilon^2{\rm d}u_t+(u_t-\Delta u){\rm d}t=f(u){\rm
d}t+\delta\phi{\rm d}W
\end{equation}
with
\[
u(0,x)=u_0(x),~u_t(0,x)=u_1(x),~u(t,0)=u(t,\pi)=0,
\]
where $\phi\in H_0^1(0,\pi)$, $u_t:=\dfrac{{\rm d}u}{{\rm d}t}$. We
assume that the nonlinear term $f$ is globally Lipschitz continuous
on $L^2(0,\pi)$ with Lipschitz constant ${\rm Lip}f$.

Rewrite (\ref{wa}) as
\[
\left\{
  \begin{array}{l}
    u_t=v, \\
    \epsilon^2v_t+v+\tilde Au=f(u)+\delta\phi\dfrac{{\rm d}W}{{\rm d}t},
  \end{array}
\right.
\]
where $\tilde Au:=-\Delta u$, and $(u,v)\in E:=H_0^1(0,\pi)\times
L^2(0,\pi)$. Let $\bar u=u$, $\bar v=v-\delta\phi z$. Here $z$
satisfies
\begin{equation}\label{zp}
\epsilon^2{\rm d}z+z{\rm d}t={\rm d}W.
\end{equation}

Let $U=(\bar u,\bar v)\in E$, then $U$ satisfies
\begin{equation}\label{V}
\dot{U}=AU+F(\theta_t\omega,U),
\end{equation}
where
\begin{equation}\label{af}
A:=\left(
           \begin{array}{cc}
             0 & {\rm id}_{L^2} \\
             -\epsilon^{-2}\tilde A & -\epsilon^{-2}{\rm id}_{L^2} \\
           \end{array}
         \right),~~
F(\omega,U):=\left(
                           \begin{array}{c}
                            \delta\phi z \\
                          \epsilon^{-2}f(\bar u)\\
                           \end{array}
                         \right).
\end{equation}
Noting that (\ref{V}) is a particular form of (\ref{main}) with
$z=0$. It is easy to verify that $A$ is the infinitesimal generator
of a $C^0$-semigroup $e^{At}$ on Hilbert space $E$. Since $F$ is
Lipschitz continuous with respect to $U$ (see (\ref{lip})), by the
classical semigroup theory concerning the local existence and
uniqueness of the solutions of evolution differential equations in
\cite{Paz}, we obtain the existence and uniqueness of (\ref{V}) and
hence (\ref{wa}).

Since the eigenvalues of $\tilde A$ are $\tilde\lambda_k=k^2$ with
corresponding eigenvectors $\tilde e_k=\sin kx$, $k=1,2,\cdots$, the
eigenvalues of the operator $A$ are
\[
\lambda_k^\pm=\frac{-1\pm\sqrt{1-4\epsilon^2k^2}}{2\epsilon^2}
\]
with corresponding eigenvectors
\[
e_k^\pm=\left(
          \begin{array}{c}
            \sin kx\\
            \lambda_k^\pm\sin kx\\
          \end{array}
        \right),~k=1,2,\cdots.
\]
It is clear that
\begin{equation}\label{lambda}
\lambda_k^+\rightarrow -k^2 ~{\rm as}~\epsilon\rightarrow 0.
\end{equation}

Denote
\[
E_1:={\rm Span}\{e_k^+|1\le k\le N\}, ~E_{-1}:={\rm
Span}\{e_k^-|1\le k\le N\},
\]
\[
E_{11}:=E_1\oplus E_{-1},~E_{22}:={\rm Span}\{e_k^\pm|k\ge N+1\},
~E_2=E_{-1}\oplus E_{22}.
\]
By the orthogonality of $\sin kx$, we have
\[
E_1\bot E_{22},~E_{-1}\bot E_{22},
\]
while $E_1$ is not orthogonal to $E_{-1}$.

Following \cite{Mor1}, we define an equivalent new inner product on
$E$. In this section, we use $\langle\cdot,\cdot\rangle$,
$\|\cdot\|$ to denote the usual inner product and norm on
$L^2(0,\pi)$, respectively. Let $U_1=(u_1,v_1)$, $U_2=(u_2,v_2)$ are
two vectors in $E$ or $E_{11}$, $E_{22}$. Recalling that the usual
inner product on $E$ defined by
\[
\langle U_1,U_2\rangle=\langle u_1,u_2\rangle+\langle \tilde
A^{\frac12}u_1,\tilde A^{\frac12}u_2\rangle+\langle v_1,v_2\rangle.
\]

Assume $\dfrac1{2\epsilon}>N+1$, define the new inner product as
follows:
\begin{align*}
\langle U_1,U_2\rangle_{E_{11}}:=& \frac1{4\epsilon^2}\langle
u_1,u_2\rangle-\langle\tilde A^{\frac12}u_1,\tilde
A^{\frac12}u_2\rangle+\langle\frac1{2\epsilon}u_1+\epsilon
v_1,\frac1{2\epsilon}u_2+\epsilon v_2\rangle,\\
\langle U_1,U_2\rangle_{E_{22}}:=& \langle\tilde
A^{\frac12}u_1,\tilde
A^{\frac12}u_2\rangle+\big(\frac1{4\epsilon^2}-2(N+1)^2\big)\langle
u_1,u_2\rangle+\langle\frac1{2\epsilon}u_1+\epsilon
v_1,\frac1{2\epsilon}u_2+\epsilon v_2\rangle.
\end{align*}
For $U=U_{11}+U_{22}$, $V=V_{11}+V_{22}$, define
\[
\langle U,V\rangle_{E}:=\langle
U_{11},V_{11}\rangle_{E_{11}}+\langle U_{22},V_{22}\rangle_{E_{22}}.
\]

Since $\dfrac1{2\epsilon}>N+1$, it is clear that $\langle
\cdot,\cdot\rangle_{E_{11}}$ is equivalent to the usual inner
product on $E_{11}$, and $\langle\cdot,\cdot\rangle_{E_{22}}$ is
equivalent to the usual inner product on $E_{22}$. Hence the new
inner product $\langle\cdot,\cdot\rangle_{E}$ is equivalent to the
usual product on $E$, see \cite{Mor1} for details.

By the definition of new inner product, it is clear that for
$U=(u,v)$ with $u=0$ we have
\begin{equation}\label{inp}
\|U\|_E=\epsilon\|v\|,
\end{equation}
and for any $U=(u,v)\in E$ we have
\begin{equation}\label{inpr}
\|U\|_E\ge\sqrt{\dfrac1{4\epsilon^2}-(N+1)^2}~\|u\|.
\end{equation}

Under the new inner product $\langle\cdot,\cdot\rangle_E$, by the
orthogonality of $\sin kx$ it is easy to verify that we have
\[
E_1\bot E_{22},~E_{-1}\bot E_{22}.
\]
Moreover, we have $E_1\bot E_{-1}$ and hence $E_1\bot E_2$. In fact,
by the definition of $\langle\cdot,\cdot\rangle_{E_{11}}$ it follows
that
\[
\left\{
  \begin{array}{l}
    \langle e_k^+,e_l^-\rangle_{E_{11}}=0, ~{\rm when}~1\le k,l\le N,~ k\neq l,\\
     \langle
e_k^+,e_k^-\rangle_{E_{11}}=\dfrac1{4\epsilon^2}-k^2
+(\dfrac1{2\epsilon}+\epsilon\lambda_k^+)(\dfrac1{2\epsilon}+\epsilon\lambda_k^-)=0,~{\rm
for}~1\le k\le N,
  \end{array}
\right.
\]
which verifies $E_1\bot E_{-1}$.

We use $A_1$, $A_2$, $A_{-1}$, $A_{22}$ to denote $A|_{E_1}$,
$A|_{E_2}$, $A|_{E_{-1}}$, $A|_{E_{22}}$, respectively. Then similar
to \cite{Mor1}, we have
\begin{align}
&\|e^{A_1t}\|=e^{\lambda_N^+t}, ~{\rm for}~t\le 0,\label{ex1}\\
&\|e^{A_{-1}t}\|=e^{\lambda_N^-t}, ~{\rm for}~t\ge 0,\label{ex2}\\
&\|e^{A_{22}t}\|=e^{\lambda_{N+1}^+t}, ~{\rm for}~t\ge 0,\label{ex3}
\end{align}
where $\|\cdot\|$ denotes the operator norm in Hilbert space
$(E,\langle\cdot,\cdot,\rangle_E)$. By (\ref{ex2}), (\ref{ex3}) we
have
\begin{equation}
\|e^{A_2t}\|=e^{\lambda_{N+1}^+t}, ~{\rm for}~t\ge 0.\label{ex4}
\end{equation}

Next we show that $F$ is Lipschitz with respect to $U$ under the
norm $\|\cdot\|_E$ and the Lipschitz constant is independent of
$\epsilon$ when $\epsilon$ is small. In fact,
\begin{align}
\|F(\omega,U_1)-F(\omega,U_2)\|_E & \le\left\|\epsilon^{-2}\left(
                                        \begin{array}{c}
                                          0\\
                                          f(u_1)-f(u_2)\nonumber\\
                                        \end{array}
                                      \right)
\right\|_E\nonumber\\
& \le^{(\ref{inp})}\epsilon^{-1}\|f(u_1)-f(u_2)\|\nonumber\\
& \le \epsilon^{-1}{\rm Lip}f\|u_1-u_2\|\nonumber\\
& \le^{(\ref{inpr})}\epsilon^{-1} \dfrac{{\rm
Lip}f}{\sqrt{\dfrac1{4\epsilon^2}-(N+1)^2}}\|U_1-U_2\|_E\nonumber\\
&\le \dfrac{{\rm
Lip}f}{\sqrt{\dfrac14-\epsilon^2(N+1)^2}}\|U_1-U_2\|_E\nonumber\\
&\le 3{\rm Lip}f\|U_1-U_2\|_E,\label{lip}
\end{align}
where the last ``=" holds when $\epsilon$ is appropriately small.

\begin{theorem}
Consider stochastic wave equation (\ref{wa}). There exists some
$\epsilon_0>0$ such that for any $\epsilon\in(0,\epsilon_0)$, the
equation (\ref{wa}) has a stochastic IM.
\end{theorem}
{\bf Proof.} Consider (\ref{V}) and let
$H=(E,\langle\cdot,\cdot\rangle_E)$, $A$ be as in (\ref{af}),
$\alpha=\lambda_N^+$, $\beta=\lambda_{N+1}^+$ and
$\eta=\dfrac{\alpha+\beta}{2}$. By (\ref{ex1}) and (\ref{ex4}), the
pseudo exponent dichotomy condition (\ref{expo}) holds with
$PH=E_1$, $QH=E_2$ and $K=1$. According to (\ref{lambda}) and
(\ref{lip}), there exists $\epsilon_0>0$ such that for any
$\epsilon\in(0,\epsilon_0)$ the spectral gap condition (\ref{gap1})
holds when $N$ is appropriately large. Hence Theorem \ref{th1} holds
for (\ref{V}), i.e. there exists a stochastic IM $M(\omega)$  for
(\ref{V}).

For $U\in E$, define the transform
\[
T(\omega,U)=U-(0,\delta\phi z),~T^{-1}(\omega,U)=U+(0,\delta\phi z).
\]
If $(t,\omega,U_0)\rightarrow \varphi(t,\omega,U_0)$ is the RDS
generated by (\ref{V}), then it is easy to verify that
\[
(t,\omega,U_0)\rightarrow\tilde\varphi:=T^{-1}(\theta_t\omega,\varphi(t,\omega,T(\omega,U_0)))
\]
is the RDS generated  by (\ref{wa}).

 Let
\[
\tilde
M(\omega):=T^{-1}(\omega,M(\omega))=\{\xi+h(\xi,\omega)+(0,\delta\phi
z)|\xi\in PE\},
\]
then $\tilde M(\omega)$ is a stochastic IM for (\ref{wa}). In fact,
\begin{align*}
\tilde\varphi(t,\omega,\tilde
M(\omega))&=T^{-1}(\theta_t\omega,\varphi(t,\omega,T(\omega,\tilde
M(\omega)))) \\
&=T^{-1}(\theta_t\omega,\varphi(t,\omega,M(\omega)))\\
&\subset
 T^{-1}(\theta_t\omega,M(\theta_t\omega))=\tilde M(\theta_t\omega),
\end{align*}
i.e. $\tilde M(\omega)$ is an invariant manifold for (\ref{wa}).

Assume $\tilde U_1$ is a solution of (\ref{wa}), then it is easy to
verify that
\[
U_1:=T(\theta_t\omega,\tilde U_1(t,\omega,T^{-1}(\omega,\tilde
U_1(0))))
\]
is a solution of (\ref{V}). By the asymptotic complete
property of $M(\omega)$, there exists a solution $U_2$ of (\ref{V})
lying on $M(\omega)$ such that
\[
\|\varphi(t,\omega,U_1(0))-\varphi(t,\omega,U_2(0))\|_E\le
c(\omega)e^{\eta t}\|U_1(0)-U_2(0)\|_E,~\forall t\ge 0.
\]
Let $\tilde U_2:=T^{-1}(\theta_t\omega,U_2(t,\omega,U_2(0)))$, then
it is easy to verify that $\tilde U_2$ is a solution of (\ref{wa})
and $\tilde U_2$ lies on $\tilde M(\omega)$. Furthermore,
\begin{align*}
\|\tilde\varphi(t,\omega,\tilde
U_1(0))-\tilde\varphi(t,\omega,\tilde
U_2(0))\|_E&=\|\varphi(t,\omega,U_1(0))-\varphi(t,\omega,U_2(0))\|_E\\
&\le  c(\omega)e^{\eta t}\|U_1(0)-U_2(0)\|_E\\
& \le c(\omega)e^{\eta t}\|\tilde U_1(0)-\tilde U_2(0)\|_E,~\forall
t\ge 0.
\end{align*}
Therefore, $\tilde M(\omega)$ has asymptotic completeness property
and hence it is a stochastic IM for (\ref{wa}). The proof is
complete.  \hfill$\Box$

\begin{remark}\rm
In above theorem we obtain the existence of stochastic IM when
$\epsilon$ is small. In fact, when $\epsilon$ is large,
counterexample has shown that the attractor of (\ref{wa}) in the
deterministic case (i.e. $\delta=0$) is not contained in any finite
dimensional manifold, see \cite{Mor3} for details. It seems that the
corresponding result holds for stochastic case, i.e. we would not
obtain the existence of stochastic IM for (\ref{wa}) when $\epsilon$
is large.
\end{remark}

Denote
\[
C_\eta^-:=\{\phi:(-\infty,0]\rightarrow E|\phi~{\rm
continuous},~\sup_{t\le 0}e^{-\eta t}\|\phi(t)\|_E<\infty\},
\]
then $C_\eta^-$ is a Banach space with norm
$\|\phi\|_{E,C_\eta^-}:=\sup_{t\le 0}e^{-\eta t}\|\phi(t)\|_E$.
Assume $R>0$ and $M_\delta(\omega)$ is a stochastic IM of
(\ref{wa}). Let
\[
M_\delta^R(\omega):=\{\xi+h(\xi,\omega)|\xi\in PE,\|\xi\|_E\le R\},
\]
where the graph of $h$ gives the IM $M_\delta(\omega)$. The
following theorem states that the stochastic IM of (\ref{wa})
converges to its deterministic counterpart almost surely when the
intensity of noise tends to zero.

\begin{theorem}
Assume $M_\delta(\omega)$ is a stochastic IM of (\ref{wa}) and $M_0$
is the IM of (\ref{wa}) when $\delta=0$ with the same dimension as
that of $M_\delta(\omega)$, then, for any $R>0$, we have
\[
\lim_{\delta\rightarrow 0}\sup_{U\in M_\delta^R(\omega)}\inf_{V_\in
M_0}\|U-V\|_E=0
\]
almost surely.
\end{theorem}
{\bf Proof.} Assume $\bar u$, $u$ satisfy
\[
\epsilon^2{\rm d}\bar u_t+(\bar u_t+\tilde A\bar u){\rm d}t=f(\bar
u){\rm d}t+\delta\phi{\rm d}W
\]
and
\[
\epsilon^2{\rm d}u_t+(u_t+\tilde A u){\rm d}t=f(u){\rm d}t,
\]
respectively. We also assume that $(\bar u,\bar u_t)$, $(u,u_t)$ lie
on $M_\delta(\omega)$, $M_0$, respectively. Let $w=\bar u-u$, then
$w$ satisfies
\begin{equation}\label{wor}
\epsilon^2w_{tt}+w_t+\tilde Aw=f(u+w)-f(u)+\delta\phi\frac{{\rm
d}W}{{\rm d}t}.
\end{equation}
Let $\bar W=(w,w_t-\delta\phi z)$, where $z$ satisfies
$\epsilon^2{\rm d}z+z{\rm d}t={\rm d}W$, then $\bar W$ satisfies
\begin{equation}\label{W}
\dot{\bar W}=A\bar W+F(\theta_t\omega,\bar W),
\end{equation}
where
\[
A=\left(
           \begin{array}{cc}
             0 & {\rm id}_{L^2} \\
             -\epsilon^{-2}\tilde A & -\epsilon^{-2}{\rm id}_{L^2} \\
           \end{array}
         \right),~~  F=\left(\begin{array}{c}
                                      \delta\phi z \\
                                       \epsilon^{-2}[f(u+w)-f(u)]
                                     \end{array}\right).
\]
It is clear that the form of (\ref{W}) is the same as that of
(\ref{V}) except that the nonlinear term $F$ is not the same. But it
is easy to verify that the nonlinear term $F$ in (\ref{W}) is
globally Lipschitz continuous with respect to $\bar W$, so (\ref{W})
has a stochastic IM and by similar argument to that of Theorem
\ref{th1} (see also (27) in \cite{Dua1}) we have $\bar W\in
C_\eta^-$ and $\bar W$ satisfies
\begin{align*}
\bar W(t)=& e^{At}P\bar
W(0)+\int_0^te^{A(t-s)}PF(\theta_s\omega,\bar W(s)){\rm
d}s\\
&\quad +\int_{-\infty}^te^{A(t-s)}QF(\theta_s\omega,\bar W(s)){\rm
d}s.
\end{align*}
Since $PE=E_1$ is of finite dimension, we can choose $(u(0),
u_t(0))$ such that $P\bar W(0)=P(\bar u(0)-u(0),\bar
u_t(0)-u_t(0)-\delta\phi z)=0$. Therefore,
\begin{align*}
e^{-\eta t}\|\bar W(t)\|_E
& \le e^{-\eta
t}\int_t^0e^{\lambda_N^+(t-s)}\left\|\left(
                                  \begin{array}{c}
                                    \delta\phi z(\theta_s\omega) \\
                                    \epsilon^{-2}[f(u+w)-f(u)]\\
                                  \end{array}
                                \right)
\right\|_E{\rm d}s\\
&\qquad +e^{-\eta
t}\int_{-\infty}^te^{\lambda_{N+1}^+(t-s)}\left\|\left(
                                  \begin{array}{c}
                                    \delta\phi z(\theta_s\omega) \\
                                    \epsilon^{-2}[f(u+w)-f(u)]\\
                                  \end{array}
                                \right)
\right\|_E{\rm d}s \\
& \le e^{-\eta t}\int_t^0e^{\lambda_N^+(t-s)}\left\|\left(
                                  \begin{array}{c}
                                    0\\
                                    \epsilon^{-2}[f(u+w)-f(u)]\\
                                  \end{array}
                                \right)
\right\|_E{\rm d}s\\
&\qquad +e^{-\eta t}\int_{t}^0e^{\lambda_{N}^+(t-s)}\left\|\left(
                                  \begin{array}{c}
                                    \delta\phi z(\theta_s\omega) \\
                                   0\\
                                  \end{array}
                                \right)
\right\|_E{\rm d}s\\
&\qquad +e^{-\eta
t}\int_{-\infty}^te^{\lambda_{N+1}^+(t-s)}\left\|\left(
                                  \begin{array}{c}
                                    0\\
                                    \epsilon^{-2}[f(u+w)-f(u)]\\
                                  \end{array}
                                \right)
\right\|_E{\rm d}s   \\
&\qquad +e^{-\eta
t}\int_{-\infty}^te^{\lambda_{N+1}^+(t-s)}\left\|\left(
                                  \begin{array}{c}
                                    \delta\phi z(\theta_s\omega) \\
                                   0\\
                                  \end{array}
                                \right)
\right\|_E{\rm d}s\\
&\le 3{\rm Lip}f\int_t^0e^{(\lambda_N^+-\eta)(t-s)}\left\| \left(
  \begin{array}{c}
    w \\
  0\\
  \end{array}
\right) \right\|_{E,C_\eta^-}{\rm d}s\\
&\qquad +3{\rm
Lip}f\int_{-\infty}^te^{(\lambda_{N+1}^+-\eta)(t-s)}\left\| \left(
  \begin{array}{c}
    w \\
  0\\
  \end{array}
\right) \right\|_{E,C_\eta^-}{\rm d}s\\
&\qquad +c_1(\omega)\left(\int_t^0e^{(\lambda_N^+-\eta)(t-s)}{\rm
d}s+\int_{-\infty}^te^{(\lambda_{N+1}^+-\eta)(t-s)}{\rm
d}s\right)\left\|\left(
              \begin{array}{c}
                \delta\phi \\
                0\\
              \end{array}
            \right)
\right\|_E\\
&\le 3{\rm Lip}f\left\| \left(
  \begin{array}{c}
    w \\
  0\\
  \end{array}
\right) \right\|_{E,C_\eta^-}\left(\dfrac{1}{\lambda_N^+-\eta}
+\dfrac{1}{\eta-\lambda_{N+1}^+}\right)\\
&\qquad +c_1(\omega)\left\|\left(
                                  \begin{array}{c}
                                    \delta\phi \\
                                   0\\
                                  \end{array}
                                \right)
\right\|_E\left(\dfrac{1}{\lambda_N^+-\eta}
+\dfrac{1}{\eta-\lambda_{N+1}^+}\right)\\
&\le 3{\rm Lip}f\|\bar
W\|_{E,C_\eta^-}\left(\dfrac{1}{\lambda_N^+-\eta}
+\dfrac{1}{\eta-\lambda_{N+1}^+}\right)\\
&\qquad +c_1(\omega)\left\|\left(
                                  \begin{array}{c}
                                    \delta\phi \\
                                   0\\
                                  \end{array}
                                \right)
\right\|_E\left(\dfrac{1}{\lambda_N^+-\eta}
+\dfrac{1}{\eta-\lambda_{N+1}^+}\right),
\end{align*}
where the third ``$\le$" holds for some $c_1(\omega)$ due to
(\ref{lip}) and the sublinear growth of $z(\theta_s\omega)$ with
respect to $s$. Hence we have
\[
\left[1-3{\rm Lip}f\left(\dfrac{1}{\lambda_N^+-\eta}
+\dfrac{1}{\eta-\lambda_{N+1}^+}\right)\right]\|\bar
W\|_{E,C_\eta^-} \le c_1(\omega)\left(\dfrac{1}{\lambda_N^+-\eta}
+\dfrac{1}{\eta-\lambda_{N+1}^+}\right)\left\|\left(
                                  \begin{array}{c}
                                    \delta\phi \\
                                   0\\
                                  \end{array}
                                \right)
\right\|_E.
\]
When $N$ is appropriately large and $\epsilon$ is appropriately
small we have
\[
3{\rm Lip}f\left(\dfrac{1}{\lambda_N^+-\eta}
+\dfrac{1}{\eta-\lambda_{N+1}^+}\right)\le\frac12,~\left(\dfrac{1}{\lambda_N^+-\eta}
+\dfrac{1}{\eta-\lambda_{N+1}^+}\right)\le 1,
\]
which implies that
\[
\|\bar W\|_{E,C_\eta^-} \le 2c_1(\omega)\left\|\left(
                                  \begin{array}{c}
                                    \delta\phi \\
                                   0\\
                                  \end{array}
                                \right)
\right\|_E.
\]
Returning back to (\ref{wor}), we let $\hat W=\bar W+(0,\delta\phi
z)=(w,w_t)$, then
\begin{align*}
\|\hat W\|_{E,C_\eta^-}& \le \|\bar W\|_{E,C_\eta^-}+\left\|\left(
                                  \begin{array}{c}
                                    0\\
                                  \delta\phi z\\
                                  \end{array}
                                \right)
\right\|_{E,C_\eta^-}\\
&\le \|\bar W\|_{E,C_\eta^-}+\left\|\left(
                                  \begin{array}{c}
                                    0\\
                                  \delta\phi\\
                                  \end{array}
                                \right)
\right\|_{E}\sup_{t\le 0}e^{-\eta t}|z(\theta_t\omega)|\\
& \le 2c_1(\omega)\left\|\left(
                                  \begin{array}{c}
                                    \delta\phi\\
                                  0\\
                                  \end{array}
                                \right)
\right\|_{E}+c_2(\omega)\left\|\left(
                                  \begin{array}{c}
                                   0 \\
                                  \delta\phi\\
                                  \end{array}
                                \right)
\right\|_{E}\\
&\le 2\delta c_1(\omega)\left\|\left(
                                  \begin{array}{c}
                                  \phi\\
                                  0\\
                                  \end{array}
                                \right)
\right\|_{E}+\delta c_2(\omega)\left\|\left(
                                  \begin{array}{c}
                                   0 \\
                                  \phi\\
                                  \end{array}
                                \right)
\right\|_{E}
\end{align*}
for some $c_2(\omega)$. Thus
\[
\|\hat W(0)\|_E\le \|\hat W\|_{E,C_\eta^-}\le\delta c_3(\omega),
\]
where
\[
c_3(\omega):=2c_1(\omega)\left\|\left(
                                  \begin{array}{c}
                                  \phi\\
                                  0\\
                                  \end{array}
                                \right)
\right\|_{E}+c_2(\omega)\left\|\left(
                                  \begin{array}{c}
                                   0 \\
                                  \phi\\
                                  \end{array}
                                \right)
\right\|_{E}.
\]
 The proof is complete. \hfill$\Box$

\section*{Acknowledgement} I am most indebted to my advisor, Professor Yong Li, not only for
his direct helpful suggestions but primarily for his continual
instruction, encouragement and support over all these years.

\section*{Appendix}

It is well-known that, for deterministic evolution equations, the
inertial manifolds contain the corresponding global attractors when
they both exist. Like deterministic case, we have the same result
for stochastic evolution equations: stochastic IM contains the
corresponding random attractor when they both exist. Here we give a
simple proof of this result.

First let us recall the definition of (global) random attractor.

\begin{definition}(\cite{Cra1})\label{attractor}
Assume $\varphi$ is an RDS on a Polish space $X$, then a random
compact set $A(\omega)$ is called a  {\em (global) random attractor}
for the RDS $\varphi$ if
\begin{itemize}
  \item $A(\omega)$ is invariant, i.e.
\begin{equation}\label{inva}
\varphi(t,\omega,A(\omega))=A(\theta_t\omega),~\forall t\ge 0
\end{equation}
for almost all $\omega\in\Omega$;

  \item $A(\omega)$ pull-back attracts every bounded deterministic
set, i.e. for any bounded deterministic set $B\subset X$, we have
\begin{equation}\label{att}
\lim_{t\rightarrow\infty}d(\varphi(t,\theta_{-t}\omega,B),A(\omega))=0
\end{equation}
almost surely.
\end{itemize}
\end{definition}
In (\ref{att}), $d(D_1,D_2)$ denotes the Hausdorff semi-metric
between $D_1$ and $D_2$, i.e.
\[
d(D_1,D_2):=\sup_{x\in D_1}\inf_{y\in D_2}d_X(x,y)
\]
for any two closed sets $D_1$, $D_2$ in $X$.

The global random attractor for RDS $\varphi$ is the {\em minimal}
random closed set which attracts all the bounded deterministic sets
and it is the {\em largest} random compact set which is invariant in
the sense of (\ref{inva}), see \cite{Cra} for details. The random
attractor defined above is unique and it is uniquely determined by
attracting {\em deterministic compact} sets, see \cite{Cr} for
details.

\begin{theorem}
Assume an SPDE has a stochastic IM $M(\omega)$ and a random
attractor $A(\omega)$. Then we have $A(\omega)\subset M(\omega)$
almost surely.
\end{theorem}
{\bf Proof.} If the assertion is false, then
\[
\mathbb P\{\omega|A(\omega)\not\subset M(\omega)\}>0.
\]
Let $\tilde A(\omega)=A(\omega)\cap M(\omega)$. Since $A(\omega)$ is
``minimal", there exists a deterministic compact set $D$ and
$\epsilon_1,\epsilon_2>0$ such that
\begin{equation}\label{con}
\mathbb
P\{\omega|\lim_{t\rightarrow\infty}d(\varphi(t,\theta_{-t}\omega,
D),\tilde A(\omega))\ge\epsilon_1\}=\epsilon_2>0.
\end{equation}
Since $A(\omega)$ is the random attractor, we have
\begin{equation}\label{con1}
\mathbb
P\{\omega|\lim_{t\rightarrow\infty}d(\varphi(t,\theta_{-t}\omega,
D),A(\omega))=0\}=1.
\end{equation}
On the other hand we have
\begin{align}
& \mathbb
P\{\omega|\lim_{t\rightarrow\infty}d(\varphi(t,\theta_{-t}\omega,
D),M(\omega))=0\}\nonumber\\
=& \mathbb P\{\omega|\lim_{t\rightarrow\infty}d(\varphi(t,\omega,
D),M(\theta_t\omega))=0\}=1\label{con2}
\end{align}
by the measure preserving of $\{\theta_t\}_{t\in\mathbb R}$ and the
fact that $M(\omega)$ is a stochastic IM for $\varphi$. According to
(\ref{con1}), (\ref{con2}) and the definition of Hausdorff
semi-metric, we have
\[
\mathbb
P\{\omega|\lim_{t\rightarrow\infty}\varphi(t,\theta_{-t}\omega,D)\subset
A(\omega)\cap M(\omega)=\tilde A(\omega)\}=1,
\]
a contradiction to (\ref{con}). The proof is complete. \hfill$\Box$

\end{document}